\documentclass[preprint]{elsarticle}
\usepackage{amsmath,amssymb,enumerate,amscd}
\input amssym.def
\usepackage[utf8]{inputenx}
\newcommand{\CP}[1]{\makebox{${\rm C}_p({#1})$}}
\newcommand{\CPB}[1]{\makebox{${\rm C}_p^{\star}({#1})$}}

\newcommand{\USC}[1]{\makebox{${\rm USC}_p({#1})$}}

\newcommand{\USCB}[1]{\makebox{${\rm USC}_p^{\star}({#1})$}}

\newcommand{\seq}[2]{\left\langle #1:#2\in\omega\right\rangle}
\newcommand{\seqc}[1]{\makebox{$[#1]_{\scriptstyle seq}$}}

\newcommand{\nula}{\makebox{${\mathbf 0}$}}
\newcommand{\constant}[1]{\makebox{${\mathbf #1}$}}
\newcommand{\snula}{\makebox{${\scriptstyle\mathbf 0}$}}
\newcommand{\RR}{{\mathbb R}}
\newcommand{\QQ}{{\mathbb Q}}
\newcommand{\emm}{\bf}
\newcommand{\pf}{{\parindent0cm{\bf Proof:\ }}}

\newtheorem{theorem}{Theorem}[section]
\newtheorem{corollary}[theorem]{Corollary}
\newtheorem{lemma}[theorem]{Lemma}
\newtheorem{problem}[theorem]{Problem}
\newtheorem{remark}[theorem]{Note}
\begin{document}
\begin{frontmatter}
\title{Selectors for dense subsets of function spaces}
\tnotetext[t1]{The work of the first author on this paper was supported by the Grant No.~1/0097/16 of Slovak Grant Agency VEGA.}
\author{Lev Bukovsk\'y}
\ead{lev.bukovsky@upjs.sk}
\author{Alexander V. Osipov}
\ead{OAB@list.ru} 
\address{Institute of Mathematics, Faculty of Science, 
P.J. \v Saf\'arik University,
Jesenn\'a~5, 040~01~Ko\v sice, Slovakia\\
Krasovskii Institute of Mathematics and Mechanics, Ural Federal University, and
 Ural State University of Economics, Yekaterinburg, Russia}
\begin{keyword}
Upper semicontinuous function\sep dense subset \sep sequentially dense subset \sep upper dense set \sep upper sequentially dense set \sep  pointwise dense subset
\sep covering propery {\rm S}${}_1$ \sep selection principle {\rm S}${}_1$.
\MSC[2010] 54C35 \sep 54C20 \sep 54D55.
\end{keyword}
\begin{abstract}
Let $\USCB{X}$ be the topological space of real upper semicontinuous bounded functions defined on $X$ with the subspace topology of the product topology on ${}^X\RR$. $\tilde\Phi^{\uparrow},\tilde\Psi^{\uparrow}$ are the sets of all upper sequentially dense, upper dense or pointwise dense subsets of $\USCB{X}$, respectively. We prove several equivalent assertions to that that $\USCB{X}$ satisfies the selection principles {\rm S}${}_1(\tilde\Phi^{\uparrow},\tilde\Psi^{\uparrow})$, including a~condition on the topological space~$X$.
\newline
We prove similar results for the topological space $\CPB{X}$ of continuous bounded functions.
\newline
Similar results hold true for the selection principles {\rm S}${}_{\scriptstyle fin}(\tilde\Phi^{\uparrow},\tilde\Psi^{\uparrow})$.
\end{abstract}
\end{frontmatter}
\section{Introduction}
We shall study the relationship between selection properties of
covers of a~topological space $X$ and selection properties of
dense subsets of the set $\USCB{X}$ of all bounded upper
semicontinuous functions on $X$ and the set of all
bounded continuous functions $C^*_p(X)$ on $X$ with the topology
of pointwise convergence. Similar problems were studied by
M.~Scheepers~\cite{Sch1,Sch4}, J.~Hale\v s~\cite{Ha}, A.~Bella, M.~Bonanzinga, M.~Matveev
~\cite{BBM}, M.~Sakai~\cite{Sa1} -- \cite{Sa5}, W.~Just, A.W.~Miller,
M.~Scheepers, P.J.~Szeptycki~\cite{JMSS}, and authors of this paper
~\cite{BL2,BL3,Osi,Osi1,Osi2,Osi3} for the set $C_p(X)$ of all
real-valued continuous functions with the topology of pointwise
convergence.
\par
We shall find equivalent conditions in the term of some properties
of density of subsets of~$\USCB{X}$ and the covering properties
S${}_1(\Phi,\Psi)$ for $\Phi,\Psi={\mathcal O},\Omega,\Gamma$. The
main results are presented in Theorems \ref{PhiD}, \ref{noAD},
\ref{PhiS}, \ref{noAS} and \ref{PhiP}. Also we present similar
Theorems \ref{PhiDcon}, \ref{PhiScon} and \ref{PhiPcon} for
bounded continuous functions.
\section{Terminology and notations}\label{TN}
\par
In the paper $( X, \tau)$ is an~infinite Hausdorff topological
space, $\tau$ is the topology, i.e., the set of open sets.
\par
We shall follow the standard terminology and notation as that of~\cite{Eng}. We recall some notions.
\par
If $A\subseteq X$ then the {\emm sequential closure $\seqc{A}$ of} $A$ is the set of all limits of sequences from $A$. 
\par
From some technical reasons, a~cover is called $o${\emm -cover}. A~cover~${\mathcal U}$ of $X$ is an~$\omega${\emm -cover}
if $X\notin {\mathcal U}$ and for any finite set $F\subseteq X$ there exists $U\in{\mathcal U}$ that contains $F$ as a subset.
An~infinite cover ${\mathcal U}$ is a~$\gamma${\emm -cover}
if for every $x\in X$ the set $\{U\in {\mathcal U}:x\notin U\}$ is finite. By ${\mathcal O}(X)$, $\Omega(X)$, $\Gamma(X)$,
or simply ${\mathcal O}$, $\Omega$, $\Gamma$ when the topological space $X$ is understood, we denote the family of all open, open $\omega$ and open $\gamma$-covers, respectively.
\par
Let us recall the a~cover~${\mathcal V}$ is a~{\emm refinement of a~cover}~${\mathcal U}$ iff
\begin{equation}\label{ref}
(\forall V\in {\mathcal V})(\exists U\in {\mathcal U}\,V\subseteq U.
\end{equation}
\par
If $\varphi$ denotes one of the symbols $o$, $\omega$ or $\gamma$, then a~$\varphi$-cover~${\mathcal U}$ is {\emm shrinkable},
if there exists an~open~$\varphi$-cover ${\mathcal W}$ such that
\begin{equation}\label{shrink}
(\forall W\in{\mathcal W})(\exists U_W\in{\mathcal U})\,
\overline{W}\subseteq U_W.
\end{equation}
The family $\{U_W:W\in{\mathcal W}\}\subseteq {\mathcal U}$ is
a~$\varphi$-cover as well. The family of all open shrinkable $o$-,
$\omega$- or $\gamma$-covers of $X$ will be denoted by ${\mathcal
O}^{\scriptstyle sh}(X)$, $\Omega^{\scriptstyle sh}(X)$ and
$\Gamma^{\scriptstyle sh}(X)$, or simply ${\mathcal
O}^{\scriptstyle sh}$, $\Omega^{\scriptstyle sh}$ and
$\Gamma^{\scriptstyle sh}$, respectively.
\par
The set ${}^X\RR$ of all real function defined on $X$ is endowed with the product topology. Thus, a~typical neighborhood of a~function \makebox{$g\in{}^X\RR$} is the set
\begin{equation}\label{neig}
V=\{h\in{}^X\RR:\vert h(x_j)-g(x_j)\vert<\varepsilon:j=0,\dots,k\}
\end{equation}
where $\varepsilon$ is a~positive real and $x_0,\dots,x_k\in X$. A~sequence of real functions $\seq{f_n}{n}$ converges to a~real function $f$ in this topology if it converges pointwise, i.e., if $\lim_{n\to\infty}f_n(x)=f(x)$ for each $x\in X$.
\par
Similarly as in~\cite{BL2} we introduce the following properties of a~family $F$ of real functions and a~function $h\in {}^X\RR$:
\[
\begin{array}{ll}
({\mathcal O}_h)&h(x)\in\overline{\{f(x):f\in F\}}\textnormal{\ for\ every\ } x\in X.\\
(\Omega_h)&h\notin F\textnormal{\  and\ }h\in\overline{F} \ in\ the\ topology\ of\ {}^X\RR.\\
(\Gamma_h)&F\textnormal{\ is\ infinite and for\ every\ }\varepsilon>0\textnormal{\ and for every\ }x\in X\\
{}&\textnormal{the\ set\ }\{f\in F:\vert f(x)-h(x)\vert\geq \varepsilon\}\textnormal{\ is\ finite}.
\end{array}
\]
Let $H\subseteq {}^X\RR$. We set
\[\Phi_h(H)=\{F\subseteq H:F\textnormal{\ possesses\ }(\Phi_h)\land (\forall f\in F)\,(f\geq h\land f-h\in H)\}.\]
\par
For a~real $a$, we denote by $\constant{a}$ the constant function on $X$ with value~$a$. For simplicity for $g\in {}^X\RR$, instead of $g+\constant{a}$ or $g-\constant{a}$ we shall write $g+a$ or $g-a$, respectively. Similarly for $\min\{\constant{a},g\}$ or $\max\{\constant{a},g\}$. If $F\subseteq {}^X\RR$, then
\[F^+=\{f\in F:f\geq 0\},\ \ F^*=\{f\in F :f\textnormal{\ is bounded}\}.\]
\par
$\CP{X}$ or $\USC{X}$ denote the set of all real continuous or upper semicontinuous functions\footnote{A function $f:X\longrightarrow \RR$ is said to be upper semicontinuous if for every real~$a$ the set $\{x\in X:f(x)<a\}$ is open.}  defined on the~topological space $X$. Instead of $\CP{X}^*$ or $\USC{X}^*$ we write $\CPB{X}$ or $\USCB{X}$, respectively.
\par
A~set $F\subseteq H\subseteq {}^X\RR$ is {\emm sequentially dense
in} $H$ if $H\subseteq \seqc{F}$. The~set $F\subseteq H\subseteq {}^X\RR$ is {\emm countably dense in} $H$
if for every function~$f\in H$ there exists a~countable set
$G\subseteq F$ such that $f\in\overline{G}$. As~obviously, the~set
$F$ is {\emm dense in} $H$  if $H\subseteq \overline{F}$.  
Finally, the~set $F$ is {\emm pointwise dense in} $A\subseteq\RR$ if $A\subseteq
\overline{\{f(x):f\in F\}}$ for each $x\in X$ ({\bf $1$-dense set} in terminology of \cite{Osi2,Osi3}). We set
\[
\begin{array}{lll}
{}&{\mathcal S}(H)=&\!\!\!\! \{F\subseteq H:F\textnormal{\ is sequentially dense in\ }H\}, \nonumber \\
{}&{\mathcal D}(H)=&\!\!\!\! \{F\subseteq H:F\textnormal{\ is dense in\ }H\} , \nonumber \\
{}&{\mathcal P}(H)=&\!\!\!\! \{F\subseteq H:F\textnormal{\ is pointwiae dense in\ }H\} . \nonumber
\end{array}
\]
Then
\[{\mathcal S}(H)\subseteq {\mathcal D}(H)\subseteq {\mathcal P}(H).\]
\par
Evidently a~sequentially dense set is countably dense as well. By Tong Theorem, see, e.g., \cite{Eng}, if $X$ is perfectly normal topological space then every (bounded) upper semicontinuous function is a~limit of a~non-increasing sequence of (bounded) continuous functions. Thus for a~perfectly normal topological space $X$  the set $\CP{X}$ is sequentially dense in $\USC{X}$. Then the set $\CPB{X}$ is sequentially dense in $\USCB{X}$ as well.
\par
The authors suppose that the answer to the following problem is negative.
\begin{problem}
Let $X$ be a~normal topological space. If $\CPB{X}$ is sequentially dense in $\USCB{X}$ then $X$ is perfectly normal?
\end{problem}
\par
A~set $F\subseteq H\subseteq {}^X\RR$ is {\emm upper sequentially dense in\ }$H$ if for every $f\in H$ there exists a~sequence $\seq{h_n}{n}$ of elements of $F$ such that $h_n\to f$, $h_n\geq f$ and $h_n-f\in H$ for each $n\in\omega$.  A~set $F\subseteq H$ is {\emm upper dense in\ }$H$ if for every $f\in H$ the set $\{h\in F:h\geq f\land h-f\in H\}$ is dense in the set~$\{h\in H:h\geq f\}$. Evidently, every upper sequentially dense in $H$ is also upper dense in $H$.
\par
One can easily see that if a~set $F\subseteq \USCB{X}$ is upper
dense in $\USCB{X}$, then for every continuous function $f$ the
set of upper semicontinuous functions \mbox{$\{h-f: h\in F\land
h\geq f\}$} is upper dense in $\USCB{X}^+$. If the~set
$F$ is upper sequentially dense in $\USCB{X}$, then for every
\mbox{$f\in \CPB{X}$} the set $\{h-f:h\in F\land h\geq f\}$ is upper
sequentially dense in $\USCB{X}^+$. We set
\[
\begin{array}{lll}
{}&{\mathcal S}^{\uparrow}(H)=&\!\!\!\! \{F\subseteq H:F\textnormal{\ is upper sequentially dense in\ }H\}, \nonumber \\
{}&{\mathcal D}^{\uparrow}(H)=&\!\!\!\! \{F\subseteq H:F\textnormal{\ is upper dense in\ }H\} , \nonumber \\
{}&{\mathcal P}^{\uparrow}(H)=&\!\!\!\! {\mathcal P}(H)\nonumber 
\end{array}
\]
Then
\[{\mathcal S}^{\uparrow}(H)\subseteq {\mathcal D}^{\uparrow}(H)\subseteq {\mathcal P}^{\uparrow}(H).\]
\par
We introduce the following notations. If $\Phi=\Gamma$ then $\tilde\Phi={\mathcal S}$. If $\Phi=\Omega$ then $\tilde\Phi={\mathcal D}$ and if $\Phi={\mathcal O}$ then $\tilde\Phi={\mathcal P}$. Similarly for~$\tilde\Phi^{\uparrow}$.
\par
Note that by definitions we have immediately
\begin{equation}\label{PhiPhih}
(\forall F\in \tilde\Phi(H)^{\uparrow})(\forall h\in H)(\exists G\subseteq F)\,G\in\Phi_h(H).
\end{equation}
\section{Covering Properties and Selection Properties}\label{S_1}
\par
We recall some notions  introduced in~\cite{Sch1} and \cite{JMSS}. Let ${\mathcal A},{\mathcal B}\subseteq {\mathcal P}(Y)$ for some set $Y$.
\par
The {\emm covering property} S${}_1({\mathcal A},{\mathcal B})$ means the following: for any sequence of sets $\seq{{\mathcal U}_n\in{\mathcal A}}{n}$ for every $n$ there exists \mbox{an~$U_n\in {\mathcal U}_n$} such that $\{U_n:n\in\omega\}\in{\mathcal B}$.
\par
The {\emm covering property} S${}_{\scriptstyle fin}({\mathcal A},{\mathcal B})$ means the following: for any sequence of sets $\seq{{\mathcal U}_n\in{\mathcal A}}{n}$, for every $n$ there exists a~finite set ${\mathcal V}_n\subseteq {\mathcal U}_n$ such that $\bigcup_{n\in\omega}{\mathcal V}_n\in{\mathcal B}$.
\par
For $\Phi,\Psi={\mathcal O},\Omega,\Gamma$, if the covering property S${}_1(\Phi(X),\Psi(X))$ holds true, we say that $X$ is an~S${}_1(\Phi,\Psi)$-space. Similarly for S${}_{\scriptstyle fin}$.
\par
If ${\mathcal F},{\mathcal G}\subseteq {\mathcal P}({}^X\RR)$ for some $X$ are sets of sets of functions, we prefer to speak about the~{\emm selection principle} S$_1({\mathcal F},{\mathcal G})$ or S${}_{\scriptstyle fin}({\mathcal F},{\mathcal G})$.
For $\Phi,\Psi={\mathcal O},\Omega,\Gamma$ and a~family $H\subseteq {}^X\RR$ we say that $H$ satisfies the selection principle S${}_1(\Phi_h,\Psi_h)$ if S${}_1(\Phi_h(H),\Psi_h(H))$ holds true. 
\par
In \cite{BL2} the author (Theorems 6.2 and 6.5 part 5))\footnote{The Theorem was formulated and proved for $\USC{X}^+$. One can easily see that it holds equally for $\USCB{X}^+$.}  has proved
\begin{theorem}[L. Bukovsk\' y]\label{Buk1}
Assume that $\Phi$ is one of the symbols $\Omega$ and $\Gamma$, and $\Psi$~is one of the symbols ${\mathcal O}$, $\Omega$, $\Gamma$.
Then for any couple $\langle \Phi,\Psi\rangle$ different from $\langle\Omega,{\mathcal O}\rangle$,
a~topological space $X$ is an~{\rm S}$_1(\Phi,\Psi)$-space if and only if $\USCB{X}$ satisfies the selection principle S${}_1(\Phi_{\snula},\Psi_{\snula})$.
\par
If a~topological space $X$ possesses the property\footnote{The property $(\varepsilon)$ was introduced and investigated in \cite{GN}: any $\omega$-cover contains a~countable $\omega$-subcover.} $(\varepsilon)$ then the equivalence holds true also for the couple $\langle\Omega,{\mathcal O}\rangle$.
\end{theorem}
\begin{corollary}[A.V. Osipov]\label{Osip1}
Assume that $\Phi$ is one of the symbols $\Omega$ and $\Gamma$, and $\Psi$~is one of the symbols ${\mathcal O}$, $\Omega$, $\Gamma$. Then for any couple $\langle \Phi,\Psi\rangle$ different from $\langle\Omega,{\mathcal O}\rangle$, a~topological space $X$ is an~{\rm S}$_1(\Phi,\Psi)$-space if and only if for every $h\in\USCB{X}$ the family 
$\USCB{X}$ satisfies the selection principle S${}_1(\Phi_h,\Psi_h)$.
\par
If a~topological space $X$ possesses the property $(\varepsilon)$ then the equivalence holds true also for the couple $\langle\Omega,{\mathcal O}\rangle$.
\end{corollary}
\pf
\par
The implication from right to left follows by Theorem~\ref{Buk1}.
\par
Let $h\in\USCB{X}$. If $F_n\in \Phi_h(\USCB{X})$ then 
\[F_n-h=\{g\in\USCB{X}^+:(\exists f\in F_n)\,g=f-h\}\]
belongs to $\Phi_{\snula}(\USCB{X})$. Thus we can apply the selection principle S${}_1(\Phi_{\snula},\Psi_{\snula})$. For every $ n\in\omega$ we obtain $g_n\in F_n-h$ such that $\{g_n:n\in\omega\}\in \Psi_{\snula}$. Then $\{g_n+h:n\in\omega\}\in\Psi_h$.      
\qed
\par
By Theorem 4.1, part 2)  of \cite{BL3}  similar results hold true also for the principle S${}_{\scriptstyle fin}$.
\begin{theorem}[L. Bukovsk\' y]\label{Buk2}
Assume that $\Phi$ is one of the symbols $\Omega$ and $\Gamma$, and $\Psi$~is one of the symbols ${\mathcal O}$, $\Omega$, $\Gamma$.
Then for any couple $\langle \Phi,\Psi\rangle$ different from $\langle\Omega,{\mathcal O}\rangle$, a~topological space $X$
is an~{\rm S}${}_{\scriptstyle fin}(\Phi,\Psi)$-space if and only if $\USCB{X}$ satisfies the selection principle
{\rm S}${}_{\scriptstyle fin}(\Phi_{\snula},\Psi_{\snula})$.
\par
If a~topological space $X$ possesses the property $(\varepsilon)$ then the equivalence holds true also for the couple $\langle\Omega,{\mathcal O}\rangle$.
\end{theorem}
\begin{corollary}\label{Osip2}
Assume that $\Phi$ is one of the symbols $\Omega$ and $\Gamma$,
and $\Psi$~is one of the symbols ${\mathcal O}$, $\Omega$,
$\Gamma$. Then for any couple $\langle \Phi,\Psi\rangle$ different
from $\langle\Omega,{\mathcal O}\rangle$, a~topological space $X$
is an~{\rm S}$_{\scriptstyle fin}(\Phi,\Psi)$-space if and only if
for every $h\in\USCB{X}$ the family $\USCB{X}$ satisfies the
selection principle S${}_{\scriptstyle
fin}(\Phi_h,\Psi_h)$.
\par
If a~topological space $X$ possesses the property $(\varepsilon)$ then the equivalence holds true also for the couple $\langle\Omega,{\mathcal O}\rangle$.
\end{corollary}
\par
We shall need also Theorem~6.3 of \cite{BL2} which reads as follows:
\begin{theorem}[L. Bukovsk\'y]\label{Th63}
Assume that $\Phi$ is one of the symbols $\Omega$, $\Gamma$ and $\Psi$ is one of the symbols ${\mathcal O}$, $\Omega$, $\Gamma$. Then for any couple $\langle\Phi,\Psi\rangle$ different from $\langle \Omega,{\mathcal O}\rangle$ a~normal topological space $X$ is an~{\rm S}$_1(\Phi^{\scriptstyle sh},\Psi)$-space if and only if $\CPB{X}$ satisfies the selection principle S${}_1(\Phi_{\snula},\Psi_{\snula})$.
\par
If a~topological space $X$ possesses the property $(\varepsilon)$ then the equivalence holds true also for the couple $\langle\Omega,{\mathcal O}\rangle$.
\end{theorem}
\par
Note the following:
\begin{eqnarray}\label{CPh}
{}&{if\ }\CPB{X}{\ satisfies\ the\ selection\ principle\ }{\rm S}_1(\Phi_{\snula},\Psi_{\snula}),\nonumber\\
{}&{\ then\ it}{\ satisfies\ }{\rm S}_1(\Phi_h,\Psi_h){\ for\ each\ }h\in\CPB{X}.
\end{eqnarray}
\section{Separability of spaces of real functions}\label{Sep}
Since the selection principle usually produces a~countable accordingly dense subset of the considered space of functions we present some conditions on the topological space for corresponding separability of the functions spaces. 
\par
We recall that the $i$-{\emm weight} $iw(X)$ of a~topological space $X$ is the
smallest infinite cardinal number $\kappa$ such that $X$ can be
mapped by a one-to-one continuous mapping onto a Tychonoff space
of the weight not greater than $\kappa$. 
\par
By the classical result of P.~Urysohn~\cite{Ur} one can easily see that a~Tychonoff space $\langle X,\tau\rangle$ has $iw(X)=\aleph_0$ if and only if there exists a~metric $\rho$ on $X$ such that $\langle X, \rho\rangle$ is separable and the topology~$\tau_{\rho}$ induced by the metric $\rho$ is weaker than the the topology $\tau$, i.e. $\tau_{\rho}\subseteq \tau$.
\par
The classical result is 
\begin{theorem}[N. Noble \cite{No}] \label{th31} Let $X$ be a Tychonoff space. The topological space $\CPB{X}$ is separable if and only if $iw(X)=\aleph_0$.
\end{theorem}
\par
Recall that $B_1(X)$ is a set of all first Baire class
 functions i.e., pointwise limits of continuous functions, defined
 on a space $X$. It is well-known that $\USCB{X}\subset B_1(X)$ for
 a perfectly normal space $X$.
\par
A~topological space $\langle X,\tau\rangle$ has the {\emm $V$-property}, if there exist a~metric $\rho$ on $X$ such that $\langle X,\rho\rangle$ is a~separable metric space, the topology $\tau_{\rho}$ induced by the metric~$\rho$ is weaker than the topology $\tau$ and any cozero set in the topology $\tau$  is an~{\rm F}${}_{\sigma}$ set in the topology~$\tau_{\rho}$.
\par
 In \cite{Ve} N.V. Velichko  proved the following result.
\begin{theorem} [N.V. Velichko] \label{th32}
Let $X$ be a~Tychonoff topological space. Then the following are equivalent
\begin{enumerate}
\item[{\rm a)}]  The topological space $\CPB{X}$ is sequentially separable.
\item[{\rm b)}] $X$ has the $V$-property.
\item[{\rm c)}] There exists a~countable set $S$ of bounded continuous functions such that $B^*_1(X)=\seqc{S}$.
\end{enumerate}
\end{theorem}
\par
Recall that a set is locally closed if it is an intersection of an
open and a closed set. A topological space $\langle X,\tau\rangle$ has the {\emm $OB$-property}, if
there  exists  a~metric $\rho$ on $X$ such that $\langle X,\rho\rangle$ is separable, the topology $\tau_{\rho}$ induced by the metric~$\rho$ is weaker than $\tau$ and any locally closed subset of $\langle X,\tau\rangle$ is an $F_{\sigma}$-set in the topology $\tau_{\rho}$.
\par
Note that $OB$-property implies $V$-property. 
\begin{theorem}\label{th220} 
A Tychonoff topological space $\langle X,\tau\rangle$ possesses the $OB$-property if and only if there exists a countable set $S\subseteq \CPB{X,\tau}$ such that $S$ is sequentially dense in $\USCB{X,\tau}$.
\end{theorem}
\pf\ \ \ 
\newline
 $(\to)$. Assume that $X$ has the $OB$-property, i.e. there a~metric $\rho$ on $X$ with corresponding properties. Let $f\in
 \USCB{X,\tau}$, $a<b$ being reals. Since
 \[(a,b)=\bigcup_{j\in\omega}\left ((-\infty, b)\cap
 [a+2^{-j}, +\infty)\right )\] 
and
\[f^{-1}((-\infty, b)\cap  [a+2^{-j}, +\infty))=f^{-1}((-\infty, b))\cap
f^{-1}([a+2^{-j}, +\infty))\]
is a locally closed subset of $X$ in the topology $\tau$, it follows that $f^{-1}(W)$ is an $F_{\sigma}$-set
in the topology $\tau_{\rho}$ for any open $W$ of $\mathbb{R}$. Thus by the Lebesgue -- Hausdorff Theorem, see e.g. \cite{Kr}, we obtain that $f\in B_1(X,\rho)$. Thus
\[\USCB{X,\tau}\subseteq B_1(X,\rho).\]
By Theorem \ref{th32}, there exists  a~countable set $S\subseteq \CPB{X,\rho}$ such that
$\seqc{S}=B^*_1(X,\rho)$. Since 
\[\CPB{X,\rho}\subset \CPB{X,\tau}\subset
 \USCB{X,\tau}\subset B_1(X,\rho),\]
it follows that $S\subset \CPB{X,\tau}$ is a countable
sequentially dense subset of $\USCB{X,\tau}$.
\newline
($\leftarrow$) Assume that $S=\{f_i : i\in \omega\}\subset
\CPB{X,\tau}$ is a sequentially dense subset of $\USCB{X,\tau}$. 
\par
Let $\tau'$ be the topology on $X$ with basis
\[\left\{\bigcap_{i\in A}f_i^{-1}((a_i,b_i)):A\in[\omega]^{<\omega}\land a_i,b_i\in\QQ\land a_i<b_i\textnormal{\ for\ }i\in A\right \}.\]
The topology $\tau'$ has a countable basis, is regular, therefore by the Urysohn Theorem~\cite{Ur} there exists a~metric $\rho$ on $X$ such that $\tau'=\tau_{\rho}$. Moreover, $f_i\in \CPB{X,\rho}$ for every $i\in \omega$.
\par
We claim that any locally closed set $D$ of $\langle X,\tau\rangle$  is an $F_{\sigma}$-set of $\langle X,\rho\rangle$.
Let $D=U\cap F$ where $U$ is an open set of $\langle X,\tau\rangle$ and $F$ is a
closed set of $\langle X,\tau\rangle$. Define the function $h$ as follows: $h(x)=0$ for $x\in U$ and
$h(x)=2$ for $x\in X\setminus U$. By the definition of~$h$, $U=h^{-1}(\{0\})$. Note that
$h\in USC^*_p(X,\tau)$. Hence there exists an~increasing sequence of integers $\seq{n_k}{ k}$ such that $\seq{f_{n_k}}{k}$ converges~to~$h$. It follows
that 
\[U=h^{-1}(\{0\})=\bigcup\limits_{j\in \omega}
\bigcap\limits_{k>j} f^{-1}_{n_k}([-1,1])\]
and hence $U$ is an $F_{\sigma}$-set of $\langle X,\rho\rangle$. 
\par
Similarly, we define the function $g$ as follows: $g(x)=2$ for $x\in F$ and $g(x)=0$ for $x\in X\setminus F$.
By the definition of $g$, $F=g^{-1}(\{2\})$. Note that $g\in USC^*_p(X,\tau)$
and hence there is an~increasing sequence of integers $\seq{m_k}{ k}$ such that $\seq{f_{m_k}}{k}$ converges~to~$g$. Equally as above we obtain
\[F=g^{-1}(\{2\})=\bigcup\limits_{j\in \omega} \bigcap\limits_{k>j}
f^{-1}_{m_k}([1, +\infty))\]
and hence $F$ is an~$F_{\sigma}$-set of $\langle X,\rho\rangle$ as well.
\par
It follows that $D=U\cap F$ is an $F_{\sigma}$-set of $\langle X,\rho\rangle$.
\qed
\begin{corollary} Let $X$ be a~Tychonoff topological space with the $OB$-property. Then $USC_p^*(X)$ is
sequentially separable.
\end{corollary}
\begin{problem} Is there a~topological space $X$ such that $USC_p^*(X)$ is
sequentially separable, but $X$ does not  have the $OB$-property?
\end{problem}
\par
Recall that a space is {\it perfect} if every open subset is an
$F_{\sigma}$-subset \cite{HM}. Note that a~topological space with the $OB$-property is perfect. However, the $OB$-property is stronger than to be perfect.
\begin{theorem} There is a perfect space which has not the $OB$-property.
\end{theorem}
\pf  Denote by $\mathbb{S}$ the Sorgenfrey line.
Then the space $\mathbb{S}^2$ is perfect (by Lemma 2.3 in
\cite{HM}). We claim that $\mathbb{S}^2$ has not the
$OB$-property. Consider the set $D=\{(x,-x): x\in \mathbb{S}\}$.
Then any subset $B$ of $D$ is a locally closed subset of
$\mathbb{S}^2$. Suppose that $f$ is a condensation from the space
$\mathbb{S}^2$ on a
 separable metric space $Y$, such that $f(B)$ is an $F_{\sigma}$-set
 of $Y$ for any locally closed set $B$ of $\mathbb{S}^2$. Then
 $|f(D)|=\mathfrak{c}$ and $f(D)$ is a separable metrizable subspace of $f(\mathbb{S}^2)$  such that for any subset $Q\subset f(D)$ the set $Q$ is an $F_{\sigma}$-set of $f(D)$, which is impossible (see, e.g., \cite{Kr}).
\qed
\section{Dense selectors of ${\rm USC^*_p(X)}$}\label{DS}
\par
We shall use the following families of functions. Let ${\mathcal U}$ be a~cover. We set
\begin{equation}\label{SU}
S({\mathcal U})=\{f_{U,g}+g:U\in{\mathcal U}\land g\in\USCB{X}\},
\end{equation}
where 
\begin{equation}\label{fUh}
f_{U,g}(x)=\left\{
\begin{array}{ll}
0&\mbox{\ if\ } x\in U,\\
1+\sup \vert g\vert&\mbox{\ otherwise}.
\end{array}
\right.
\end{equation}
We show the basic properties of the families $S({\mathcal U})$.
\begin{lemma}\label{lem1}\ \ \ 
\begin{enumerate}
\item[{\rm 1)}] If ${\mathcal U}$ is an~open $\omega$-cover, then the family $S({\mathcal U})\subseteq\USCB{X}$ is upper dense in $\USCB{X}$.
\item[{\rm 2)}] If ${\mathcal U}$ is an~open $\gamma$-cover, then the family $S({\mathcal U})\subseteq\USCB{X}$ is upper sequentially dense in $\USCB{X}$.
\end{enumerate}
\end{lemma}
\pf
One can easily see that $f_{U,g}+g\geq g$ and $f_{U,g}+g$ is bounded upper
semicontinuous for $g\in \USCB{X}$.\par
We show that if ${\mathcal U}$ is an~open $\omega$-cover then $S({\mathcal U})$ is upper dense in
$\USCB{X}$. Assume that $g\in\USCB{X}$. If $V$ is the neighborhood
of $g$ defined by \eqref{neig}, then there exists a~$U\in
{\mathcal U}$ such that $x_0,\dots,x_k\in U$. Then
$f_{U,g}(x_j)+g(x_j)=g(x_j)$ for $j=0,\dots,k$. Hence $f_{U,g}+g\in V$,
$f_{U,g}+g\geq g$ and $(f_{U,g}+g)-g\in\USCB{X}$. Thus $S({\mathcal U})$ is upper dense.
\par
If ${\mathcal U}$ is an~open $\gamma$-cover then $S({\mathcal U})$ is upper sequentially dense in
$\USCB{X}$. Indeed, let $g\in\USCB{X}$. Let $\{U_i:i\in\omega\}$ be a~countable $\gamma$-subcover of ${\mathcal U}$. For $i\in\omega$, we let $g_i=f_{U_i,g}+g\in S({\mathcal U})$.
We show that the sequence $\seq{g_i}{i}$ converges to $g$. Let $V$
be a~neighborhood of $g$ defined by \eqref{neig}. Since $\{U_i:i\in\omega\}$ is a~$\gamma$-cover, there exists an~$i_0$ such that $x_0,\dots,x_k\in U_i$ for $i\geq i_0$. Then for such an~$i$ we
have $g_i(x_j)=g(x_j)$ for $j=0,\dots,k$. Therefore the elements
of the sequence $\seq{g_i}{i}$ belong to $V$ for $i\geq i_0$. As above, $g_i-g\in\USCB{X}$. Thus
$S({\mathcal U})$ is upper sequentially dense.
\qed
\begin{theorem}\label{PhiD}
Let $\Phi=\Omega,\Gamma$. Then 
the following are equivalent:
\begin{enumerate}
\item[{\rm a)}] $\USCB{X}$ satisfies the selection principle {\rm
S}${}_1(\tilde \Phi^{\uparrow},{\mathcal D})$. 
\item[{\rm b)}] $\USCB{X}$ is separable and the topological space~$X$ possesses the covering property {\rm S}${}_1(\Phi,\Omega)$. 
\item[{\rm c)}] $\USCB{X}$ is separable and
satisfies the selection principle
\newline
{\rm S}${}_1(\Phi_h,\Omega_h)$ for every $h\in\USCB{X}$.
\item[{\rm d)}] $\USCB{X}$ is separable and satisfies the selection principle
\newline
{\rm S}${}_1(\tilde \Phi^{\uparrow},\Omega_h)$ for every $h\in\USCB{X}$.
\end{enumerate}
\end{theorem}
\pf
\newline
${\rm a)}\to {\rm b)}$. Let $\{{\mathcal U}_n:n\in\omega\}\subseteq \Phi$. We may assume that
${\mathcal U}_{n+1}$ is a~refinement of ${\mathcal U}_n$ for each $n\in\omega$. If $\Phi=\Gamma$
we may also assume that  for every $n\in\omega$, the cover ${\mathcal U}_n$ is a~countable family $\{U_i^n:i\in\omega\}$.
\par
 For every $n\in\omega$ we set
\begin{equation}\label{Sn}
S_n=S({\mathcal U}_n).
\end{equation}
By Lemma \ref{lem1} we have $S_n\in\tilde\Phi^{\uparrow}$. Thus, by the selection principle {\rm S}${}_1(\tilde
\Phi^{\uparrow},{\mathcal D})$, for every $n\in\omega$ we obtain
an~$f_{U_n,h_n}\in S_n$ such that $\{f_{U_n,h_n}:n\in\omega\}$
is dense in $\USCB{X}$. We show that $\{U_n:n\in\omega\}$ is
an~$\omega$-cover.
\par
Let $x_0,\dots,x_k\in X$. Consider the open non-empty set
\[U=\{g\in\USCB{X}:\vert g(x_j)\vert<1/2\textnormal{\ for\ } j=0,\dots,k\}\]
Since the set  $\{f_{U_n,h_n}:n\in\omega\}$ is dense in
$\USCB{X}$, there exists an~$n$ such that $f_{U_n,h_n}\in U$.
Since $\vert f_{U_n,h_n}(x_i)\vert <1/2$ for $i=0,\dots,k$, by \eqref{fUh} we obtain
$x_0,\dots,x_k\in U_n$.
\par 
The implication ${\rm b)}\to{\rm c)}$ follows by Corollary \ref{Osip1}.
\par
The implication ${\rm c)}\to{\rm d)}$ is obvious by \eqref{PhiPhih}.
\newline
${\rm d)}\to{\rm a)}$. We assume that $\USCB{X}$ is separable and satisfies the
selection principle {\rm S}${}_1(\tilde \Phi^{\uparrow},\Omega_h)$ for every $h\in\USCB{X}$. Thus, there exists a~ countable set $D=\{d_n:n\in\omega\}$ dense in $\USCB{X}$.
 Let $\{S_{n,m}:n,m\in\omega\}$ be a~sequence of
subsets of $\USCB{X}$ such that $S_{n,m}\in \tilde
\Phi^{\uparrow}$ for each $n,m\in \omega$. For every $n\in\omega$
we apply the sequence selection principle {\rm S}${}_1(\tilde
\Phi^{\uparrow},\Omega_{d_n})$ to the sequence $\seq{S_{n,m}}{m}$
and for every $m\in\omega$ we obtain $d_{n,m}\in S_{n,m}$ such
that $d_n\in\overline{\{d_{n,m}:m\in\omega\}}$. Then
$\{d_{n,m}:n,m\in\omega\}$ is  dense in $\USCB{X}$.
\qed
\par
By Theorem \ref{Buk1} and analogously to the proof of Theorem
\ref{PhiD} we get the following theorem.
\begin{theorem}\label{PhiD1}
Let $\Phi=\Omega$ or $\Phi=\Gamma$. Assume that $\CPB{X}$ is countably dense in $\USCB{X}$.Then for any couple $\langle\Phi,\Psi\rangle$ different from $\langle \Omega,{\mathcal O}\rangle$, the following are equivalent:
\begin{enumerate}
\item[{\rm a)}] $\USCB{X}$ satisfies the selection principle
{\rm S}${}_1(\tilde \Phi^{\uparrow},{\mathcal D})$,
\item[{\rm b)}] $\USCB{X}$ is separable  and the topological space~$X$ possesses the covering property
{\rm S}${}_1(\Phi,\Omega)$. 
\item[{\rm c)}] $\USCB{X}$ is separable and satisfies the selection principle
\newline
{\rm S}${}_1(\Phi_{\snula},\Omega_{\snula})$. 
\item[{\rm d)}] $\USCB{X}$ is separable and satisfies the selection principle
\newline
{\rm S}${}_1(\tilde \Phi^{\uparrow},\Omega_{\snula})$.
\end{enumerate}
\end{theorem}
\pf
\par
We prove only the implication ${\rm d)}\to{\rm a)}$. The proofs of other implications are almost equal to those in the proof of Theorem~\ref{PhiD}.
\par
Assume that $\CPB{X}$ is countably dense in $\USCB{X}$, $\USCB{X}$ is separable and satisfies the selection principle S${}_1(\tilde \Phi^{\uparrow},\Omega_{\snula})$. Thus, there exists a~ countable set $D=\{d_n:n\in\omega\}$ dense in $\USCB{X}$. Since $\CPB{X}$ is countably dense in $\USCB{X}$, for every $n\in\omega$ there exists a~countable set $D_n=\{d_{n,m}:m\in\omega\}\subseteq \CPB{X}$ such that $d_n\in \overline{D_n}$ for each $n\in\omega$.
\par
Let $\{S_{n,m,l}:n,m,l\in\omega\}$ be a~sequence of subsets of $\USCB{X}$, each $S_{n,m,l}$ being in $\tilde\Phi^{\uparrow}$. We can apply the sequence selection principle S${}_1(\tilde\Phi^{\uparrow},\Omega_{\snula})$ to the sequence 
\[\seq{\{h-d_{n,m}:h\in S_{n,m,l}\}}{l}.\] 
For every $l\in\omega$ we obtain $d_{n,m,l}\in S_{n,m,l}$ such that 
\[\constant{0}\in\overline{\{d_{n,m,l}-d_{n,m}:l\in\omega\}}.\]
Then
\[d_{n,m}\in\overline{\{d_{n,m,l}:l\in\omega\}}.\]
Thus  $\{d_{n,m,l}:n,m,l\in\omega\}$ is the desired  countable dense set.
\qed
\par
Since no (infinite Hausdorff) topological space has the covering property S${}_1({\mathcal O},\Omega)$, we obtain
\begin{theorem}\label{noAD}
$\USCB{X}$ does not have the property {\rm S}${}_1({\mathcal
P},{\mathcal D})$ for any topological space $X$.
\end{theorem}
\pf As in proof of Theorem~\ref{PhiD} one can show that if
$\USCB{X}$ possesses the property {\rm S}${}_1({\mathcal
P},{\mathcal D})$ then $X$ possesses the covering property
S${}_1({\mathcal O},\Omega)$. 
\qed
\section{Sequentially dense selectors of ${\rm USC^*_p(X)}$}\label{SDS}
\medskip
\begin{theorem}\label{PhiS}
Let $\Phi=\Omega$ or $\Phi=\Gamma$. Then for any couple $\langle\Phi,\Psi\rangle$ different from $\langle \Omega,{\mathcal O}\rangle$, the following are equivalent:
 \begin{enumerate}
\item[{\rm a)}] $\USCB{X}$ satisfies the selection principle {\rm
S}${}_1(\tilde \Phi^{\uparrow},{\mathcal S})$. 
\item[{\rm b)}] $\USCB{X}$ is sequentially separable and the topological space~$X$ possesses the
covering property {\rm S}${}_1(\Phi,\Gamma)$. 
 \item[{\rm c)}]  $\USCB{X}$ is sequentially separable and satisfies the selection
principle {\rm S}${}_1(\Phi_h,\Gamma_h)$ for every $h\in\USCB{X}$. 
\item[{\rm d)}] $\USCB{X}$ is sequentially separable and satisfies the selection principles 
{\rm S}${}_1(\Gamma_{\snula},\Gamma_{\snula})$ and {\rm S}${}_1(\tilde \Phi^{\uparrow},\Gamma_h)$ for every $h\in\USCB{X}$.
\end{enumerate}
\end{theorem}
\pf
\newline
${\rm a)}\to {\rm b)}$. Let $\{{\mathcal U}_n:n\in\omega\}\subseteq \Phi$. We may assume that
${\mathcal U}_{n+1}$ is a~refinement of ${\mathcal U}_n$ for each $n\in\omega$. If $\Phi=\Gamma$
we may also assume that  for every $n\in\omega$, the cover ${\mathcal U}_n=\{U_i^n:i\in\omega\}$ is a~countable family.
\par
We define the sets $S_n$ by \eqref{Sn}. By Lemma
\ref{lem1}, $S_n\in\tilde\Phi^{\uparrow}(\USCB{X})$. We apply the selection
principle {\rm S}${}_1(\tilde\Phi^{\uparrow},{\mathcal S})$ and for
every $n$ we obtain a~function~$f_{U_n,h_n}\in S_n$ such that
$\{f_{U_n,h_n}:n\in\omega\}$ is sequentially dense in
$\USCB{X}$. For every $n$ we shall find a~set $V_n\in{\mathcal U}_n$ such that  $\{V_n:n\in\omega\}$ is a~$\gamma$-cover.
\par
Evidently there exists an~increasing sequence $\seq{n_k}{k}$ such that $f_{U_{n_k},h_{n_k}}\to \nula$. We set $V_{n_k}=U_{n_k}$. If $n_k<n<n_{k+1}$, then by~\eqref{ref} we can find a~set $V_n\in{\mathcal U}_n$ such that $U_{n_{k+1}}\subseteq V_n$.
\par
Let $x\in X$. Since
\[W=\{g\in\USCB{X}:\vert g(x)\vert<1\}\]
is a~neighborhood of $\nula$, there exists an~$k_0$ such that $f_{U_{n_k},h_{n_k}}\in W$ for each $k\geq k_0$. 
If $g\in W$ then $g(x)<1$. Thus for $k\geq k_0$ we have $f_{U_{n_k},h_{n_k}}(x)<1$. Therefore $x\in U_{n_k}=V_{n_k}$. By the choose of $V_n$ for $n\notin\{n_k:k\in\omega\}$ we obtain that $x\in V_n$ for each $n\geq n_{k_0}$.
\par
The implication ${\rm b)}\to{\rm c)}$ follows by Corollary \ref{Osip1}.
\par
The implication ${\rm c)}\to{\rm d)}$ is obvious by \eqref{PhiPhih}.
\par
We prove the implication ${\rm d)}\to{\rm a)}$.
\par
Assume that there exists a~countable set
$\{d_n:n\in\omega\}\subseteq \USCB{X}$ sequentially dense in $\USCB{X}$ and $\USCB{X}$ satisfies the selection principles S${}_1(\Gamma_{\snula},\Gamma_{\snula})$ and  {\rm S}${}_1(\tilde \Phi^{\uparrow},\Gamma_h)$ for each $h\in\USCB{X}$.
\par
Let $(S_{n,m}:n,m\in\omega)$ be a~sequence of subsets of
$\USCB{X}$ all being in $\tilde\Phi^{\uparrow}$. For every $n$ we apply the selection principle {\rm S}${}_1(\tilde \Phi^{\uparrow},\Gamma_{d_n})$
to the sequence $(S_{n,m}:m\in\omega)$. Then for every $m\in\omega$ we obtain
$d_{n,m}\in S_{n,m}$, $d_{n,m}\geq d_n$, $d_{n,m}-d_n\in\USC{X}$ and such
that $d_{n,m}-d_n\to \nula$ ($m\rightarrow \infty$).
\par
We show that the set  $\{d_{n,m}:n,m\in\omega\}$ is the desired sequentially dense set.
\par
Indeed, if $h\in\USCB{X}$ then there exists an~increasing sequence $\seq{n_k}{k}$ such that $d_{n_k}\to h$.  
Since for every $k$ we have
\[d_{n_k,m}-d_{n_k}\to \nula,\]
by {\rm S}${}_1(\Gamma_{\snula},\Gamma_{\snula})$ there exists a~sequence $\seq{m_k}{k}$ such that 
\[d_{n_k,m_k} -d_{n_k}\to \nula.\]
Thus \[d_{n_k,m_k}=(d_{n_k,m_k}- d_{n_k})+d_{n_k}\]
converges to $h$ ($k\rightarrow \infty$).
\qed
\begin{corollary}
Let $\Phi=\Omega$ or $\Phi=\Gamma$. Assume that $X$ is a~Tychonoff topological space and $\CPB{X}$ is seuqentially dense on $\USCB{X}$. Then the following are equivalent:
 \begin{enumerate}
\item[{\rm a)}] $\USCB{X}$ satisfies the selection principle {\rm
S}${}_1(\tilde\Phi^{\uparrow},{\mathcal S})$. 
\item[{\rm b)}] The topological space~$X$ possesses the OB-property and the
covering property {\rm S}${}_1(\Phi,\Gamma)$. 
 \item[{\rm c)}]  $\USCB{X}$ is sequentially separable and satisfies the selection
principle {\rm S}${}_1(\Phi_h,\Gamma_h)$ for each $h\in\USCB{X}$. 
\item[{\rm d)}] $\USCB{X}$ is sequentially separable and satisfies the selection principles 
{\rm S}${}_1(\Gamma_{\snula},\Gamma_{\snula})$ and {\rm S}${}_1(\tilde \Phi^{\uparrow},\Gamma_h)$ for each $h\in\USCB{X}$.
\end{enumerate}
\end{corollary}
\pf
The only non-trivial implication is "b) of Theorem~\ref{PhiS}" to "b) of the Corollary". By S${}_1(\Phi,\Gamma)$ we obtain that $\USCB{X}$ satisfies the selection principle S${}_1(\Gamma_{\snula},\Gamma_{\snula})$. Thus we obtain a~countable subset of $\CPB{X}$ dense in $\USCB{X}$. Since $X$ is Tychonoff we can apply Theorem~\ref{th220}.
\qed  
\begin{corollary}
 Let $\Phi=\Omega$ or $\Phi=\Gamma$. Assume that the Tychonoff topological space~$X$ has the $OB$-property. Then for the following are equivalent:
 \begin{enumerate}
\item[{\rm a)}] $\USCB{X}$ satisfies the selection principle {\rm
S}${}_1(\tilde\Phi^{\uparrow},{\mathcal S})$. 
\item[{\rm b)}] $X$
possesses the covering property {\rm S}${}_1(\Phi,\Gamma)$.
\item[{\rm c)}] $\USCB{X}$ satisfies the selection principle {\rm
S}${}_1(\Phi_h,\Gamma_h)$ for every $h\in\USCB{X}$.
\item[{\rm d)}] $\USCB{X}$ satisfies the selection principles {\rm S}${}_1(\Gamma_{\snula},\Gamma_{\snula})$ and 
\newline
{\rm S}${}_1(\tilde\Phi^{\uparrow},\Gamma_{\snula})$.
\end{enumerate}
\end{corollary}
\pf
We have only to prove the implication ${\rm d)}\to{\rm a)}$. By Theorem~\ref{th220} there exists a~countable set  
$\{d_n:n\in\omega\}\subseteq \CPB{X}$ sequentially dense in $\USCB{X}$. 
\par
Let $(S_{n,m}:n,m\in\omega)$ be a~sequence of subsets of
$\USCB{X}$ all in $\tilde\Phi^{\uparrow}$. For every $n$ we apply the selection principle {\rm S}${}_1(\tilde \Phi^{\uparrow},\Gamma_{\snula})$
to the sequence $(S_{n,m}-d_n:m\in\omega)$. Then for every $m\in\omega$ we obtain
$d_{n,m}\in S_{n,m}$, $d_{n,m}\geq d_n$, $d_{n,m}-d_n\in\USC{X}$ and such
that $d_{n,m}-d_n\to \nula$ ($m\rightarrow \infty$).
\par
We can show that the set  $\{d_{n,m}:n,m\in\omega\}$ is sequentially dense subset of $\USCB{X}$ equally as in the proof of Theorem~\ref{PhiS}.
\qed
\par
Since no infinite Hausdorff topological space satisfies S${}_1({\mathcal O},\Gamma)$ as above one can easily prove
\begin{theorem}\label{noAS}
$\USCB{X}$ does not have the property {\rm S}${}_1({\mathcal
P},{\mathcal S})$ for any (infinite Hausdorff) topological space
$X$.
\end{theorem}

\section{Pointwise dense selectors of ${\rm USC^*_p(X)}$}\label{PDS}
\begin{theorem}\label{PhiP}
For $\Phi=\Gamma$ the following are equivalent:
\begin{enumerate}
\item[{\rm a)}] $\USCB{X}$ satisfies the selection principle {\rm
S}${}_1(\tilde\Phi^{\uparrow},{\mathcal P})$. 
\item[{\rm b)}] The topological space~$X$ possesses the covering property\
\newline
{\rm S}${}_1(\Phi,{\mathcal O})$.
\item[{\rm c)}] $\USCB{X}$ satisfies the selection principle
{\rm S}${}_1(\Phi_{\snula},{\mathcal O}_{\snula})$. 
\item[{\rm  d)}] $\USCB{X}$ satisfies the selection principle {\rm S}${}_1(\tilde\Phi^{\uparrow},{\mathcal O}_{\snula})$.
\end{enumerate}
If $X$ possesses the property $(\varepsilon)$ then the equivalences hold true also for $\Phi=\Omega$.
\end{theorem}
\pf
\par
${\rm a)}\to {\rm b)}$. Let $\{{\mathcal
U}_n:n\in\omega\}\in \Phi$. We can assume that for every
$n\in\omega$ the cover ${\mathcal U}_{n+1}$ is a~refinement of
${\mathcal U}_n$. If $\Phi=\Gamma$ we may also assume that  for
every $n\in\omega$, the cover ${\mathcal U}_n$ is a~countable family $\{U_i^n:i\in\omega\}$.
\par
We define the sets $S_n$ by \eqref{Sn}. By Lemma~\ref{lem1}, every $S_n$ is upper sequentially dense and pointwise
dense in $\mathbb{R}$. By the selection principle {\rm
S}${}_1(\tilde\Phi^{\uparrow},{\mathcal P})$, for every
$n\in\omega$ we obtain an~$f^n_{U_n,h_n}\in S_n$ such that the set 
$\{f^n_{U_n,h_n}:n\in\omega\}$ is pointwise dense in $\mathbb{R}$.
We show that $\{U_n:n\in\omega\}$ is a~cover.
\par
 If $x\in X$ then there exists an~$n$ such that $f^n_{U_n,h_n}(x)\in(-1/2,1/2)$. Then $x\in U_n$.
\par
By Corollary~\ref{Osip1} we obtain ${\rm b)}\to {\rm c)}$.
The implication ${\rm c)}\to {\rm d)}$ is trivial. We show the
implication ${\rm d)}\to {\rm a)}$.
\par
Let $\{S_{n,m}:n,m\in\omega\}$ be a~sequence of upper sequentially
dense subsets of $\USCB{X}$. Let the rational
numbers $\QQ$ be enumerated as $\{q_n:n\in\omega\}$. Fix $n\in\omega$. By {\rm
S}${}_1(\tilde\Phi^{\uparrow},{\mathcal O}_{\snula})$ for every
$m\in\omega$ there exists $f_{n,m}\in S_{n,m}$, $f_{n,m}\geq q_n$
such that $\{f_{n,m}-q_n: m\in\omega\}\in {\mathcal O}_{\snula}$.
Then $\{f_{n,m}:n,m\in\omega\}$ is pointwise dense in
$\mathbb{R}$.
\qed
\begin{remark}
 By Theorem \ref{Buk2} and Corollary \ref{Osip2} similar results hold true also for the principle
${\rm S}_{\scriptstyle fin}$.
\end{remark}
\section{Dense selectors of the space  $\CP{X}$}
In papers \cite{BBM,Osi,Osi2,Osi3,Sa2,Sch4}  the authors have got characterizations of dense, sequentially dense and pointwise dense
selectors of the space of real-valued continuous functions
$C_p(X)$ for Tychonoff space $X$. We can show similar results for
selectors of sequences of dense, sequentially dense and pointwise
dense subsets of $\CPB{X}$ for a~normal topological space~$X$.
\par
In the next we assume that $X$ is a~normal topological space.
\par
Let ${\mathcal U}$ be a~shrinkable cover of $X$, ${\mathcal W}$ being a~cover of $X$ such that \eqref{shrink} holds true.
We set
\begin{equation}\label{Tn}
T({\mathcal U})=\{g_{W,h}+h: W\in{\mathcal W}\land h\in \CPB{X}\},
\end{equation}
where $f_{W,h}\in\CPB{X}$ is such that
\begin{equation}\label{fW}
g_{W,h}(x)=\left\{
\begin{array}{ll}
0&\mbox{\ if\ } x\in W,\\
1+\sup \vert h\vert&\mbox{\ if\ }x\in X\setminus U_W.
\end{array}
\right.
\end{equation}
\par
The existence of such a function follows from that that $X$ is a~normal topological space.
\par
Equally as in Section~\ref{DS} one can easily prove
\begin{lemma}\label{lem5}\ \ \ 
\begin{enumerate}
\item[{\rm 1)}] If ${\mathcal U}$ is an~open shrinkable $\omega$-cover, then the family $S({\mathcal U})$ is upper dense in $\CPB{X}$.
\item[{\rm 2)}] If ${\mathcal U}$ is an~open shrinkable $\gamma$-cover, then the family $S({\mathcal U})$ is upper sequentially dense in $\CPB{X}$.
\end{enumerate}
\end{lemma}
\begin{theorem}\label{PhiDcon}
Let $\Phi=\Omega$ or $\Phi=\Gamma$. Then for a~normal topological
space~$X$ the following are equivalent:
\begin{enumerate}
\item[{\rm a)}] $\CPB{X}$ satisfies the selection principle {\rm
S}${}_1(\tilde \Phi^{\uparrow},{\mathcal D})$.
\item[{\rm b)}] $iw(X)=\aleph_0$ and the topological space~$X$ possesses the covering property {\rm
S}${}_1(\Phi^{\scriptstyle sh},\Omega)$.
\item[{\rm c)}] $\CPB{X}$ is separable and satisfies the selection principle\\
{\rm S}${}_1(\Phi_{\snula},\Omega_{\snula})$.
\item[{\rm d)}] $\CPB{X}$ is separable and satisfies the selection principle\\
 {\rm S}${}_1(\tilde \Phi^{\uparrow},\Omega_{\snula})$.
\end{enumerate}
\end{theorem}
\pf We show ${\rm a)}\to{\rm b)}$. Let $\seq{{\mathcal U}_n}{n}$ be a sequence
 of open shrinkable $\varphi$-covers of $X$.  Let $\seq{{\mathcal W}_n}{n}$ be a~sequence of open 
$\varphi$-covers such that \eqref{shrink} holds true for each $n$. We can assume that ${\mathcal U}_n=\{U_W:W\in {\mathcal W}_n\}$, where $U_W$ is the set of \eqref{shrink}.
 For every $n\in\omega$ we set
\begin{equation}\label{Tn}
T_n=T({\mathcal U}_n).
\end{equation}
By Lemma~\ref{lem5} we can apply {\rm S}${}_1(\tilde\Phi^{\uparrow},{\mathcal D})$
to the sequence $\seq{T_n}{n}$ and we obtain for every
$n\in\omega$ a~function $g_n\in T_n$ such that the set
$\{g_n:n\in\omega\}$ is dense in $\CP{X}$. By definition of
$T_n$, there exists $h_n\in\CP{X}$ and $W_n\in{\mathcal W}_n$
such that $g_n=g_{W_n,h_n}+h_n$. We show that $\{U_{W_n}:n\in\omega\}$
is an~$\omega$-cover.
\par
Let $x_0,\dots,x_k\in X$. Let us consider the neighborhood of
$\nula$ defined as \[V=\{h\in\CP{X}:\vert h(x_i)\vert<\frac 12\textnormal{\ for\ }i=0,\dots,k\}.\]
Since $\{g_n:n\in\omega\}$ is dense in $\CP{X}$ there exists
an~$n\in\omega$ such that $g_n\in V$. We know that
$g_n=g_{W_n,h_n}+h_n$. For $i=0,\dots,k$ we have $\vert g_n(x_i)\vert<\frac12$.
If $x\in X\setminus U_{W_n}$ then $g_n(x)\geq 1$. Thus $x_0,\dots,x_k\in U_{W_n}$.
\par
Similarly to the proof of Theorem~\ref{PhiD} the implication ${\rm
b)}\to{\rm c)}$ is special case of Theorem~\ref{Buk2} and the
implication ${\rm c)}\to{\rm d)}$ is trivial. We show ${\rm
d)}\to{\rm a)}$.
\par
Assume that  $\CP{X}$ is separable and satisfies the selection
principle  {\rm S}${}_1(\tilde \Phi^{\uparrow},\Omega_{\snula})$.
Let $\{d_n:n\in\omega\}$ be a~dense subset of $\CP{X}$. Let
$\{S_{n,m}:n,m\in\omega\}$ be a~sequence of subsets of
$\CP{X}$, each $S_{n,m}$ being in $\tilde\Phi^{\uparrow}$. We can
apply the sequence selection principle
S${}_1(\tilde\Phi^{\uparrow},\Omega_{\snula})$ to the sequence
\[\seq{\{h-d_n:h\in S_{n,m}\}}{m}.\]
For every $m\in\omega$ we obtain $d_{n,m}\in S_{n,m}$ such that
\[\constant{0}\in\overline{\{d_{n,m}-d_n:m\in\omega\}}.\]
Then
\[d_n\in\overline{\{d_{n,m}:m\in\omega\}}.\]
Thus  the selector $\{d_{n,m}:n,m\in\omega\}$ is a~countable dense subset of~$\CPB{X}$, every $d_{n,m}$ belongs to $S_{n,m}$.
\qed
\par
Modifying the proof of Theorem~\ref{PhiS} similarly as the proof
of Theorem~\ref{PhiDcon}, one can easily prove the next two
results. Namely
\begin{theorem}\label{PhiScon} For $\Phi={\mathcal O},\,\Omega,\,\Gamma$  and a~normal topological space $X$
the following are equivalent:
\begin{enumerate}
\item[{\rm a)}] $\CPB{X}$ satisfies the selection principle {\rm
S}${}_1(\tilde \Phi^{\uparrow},{\mathcal S})$. 
\item[{\rm b)}] The topological space~$X$ possesses the $V$-property and the covering property~{\rm S}${}_1(\Phi^{\scriptstyle sh},\Gamma)$. 
\item[{\rm c)}] $\CPB{X}$ is sequentially separable and satisfies the selection principle {\rm S}${}_1(\Phi_{\snula},\Gamma_{\snula})$. 
\item[{\rm d)}] $\CPB{X}$ is sequentially separable and satisfies
the selection principle {\rm S}${}_1(\tilde \Phi^{\uparrow},\Gamma_{\snula})$.
\end{enumerate}
\end{theorem}
and
\begin{theorem}\label{PhiPcon}
For $\Phi={\mathcal O},\,\Gamma$  and a~normal topological space
$X$ the following are equivalent:
\begin{enumerate}
\item[{\rm a)}] $\CPB{X}$ satisfies the selection principle {\rm S}${}_1(\tilde\Phi^{\uparrow},{\mathcal P})$. 
\item[{\rm b)}] $X$ possesses the covering property {\rm S}${}_1(\Phi^{\scriptstyle sh},{\mathcal O})$. 
\item[{\rm c)}] $\CPB{X}$ satisfies the selection principle {\rm S}${}_1(\Phi_{\snula},{\mathcal O}_{\snula})$. 
\item[{\rm d)}] $\CPB{X}$ satisfies the selection principle {\rm S}${}_1(\tilde\Phi^{\uparrow},{\mathcal O}_{\snula})$.
\end{enumerate}
If $X$ possesses the property $(\varepsilon)$ then the
equivalences hold true also for $\Phi=\Omega$.
\end{theorem}
\begin{remark}
Similar results hold true also for the principle
S${}_{\scriptstyle fin}$.
\end{remark}
\section{Some consequences}
Since in a~regular topological space every open
$\omega$-cover is shrinkable, see e.g. \cite{BL2}, Lemma 7.1,  in
b) of Theorems \ref{PhiDcon}, \ref{PhiScon} and \ref{PhiPcon} we can replace the corresponding covering
principle {\rm S}${}_1(\Omega^{\scriptstyle sh},\Psi)$ by {\rm S}${}_1(\Omega,\Psi)$. Hence we obtain, e.g. from
Theorems~\ref{PhiD} and \ref{PhiDcon}, the following
\begin{corollary}
For $\Psi=\Omega,\Gamma$ and a~normal topological space $X$ the following statements are
equivalent:
\begin{enumerate}
\item[{\rm a)}] $\USCB{X}$ satisfies the selection principle {\rm S}${}_1(\tilde \Omega^{\uparrow},\tilde\Psi)$. 
\item[{\rm b)}] The topological space~$X$ possesses the OB-property and the covering
property {\rm S}${}_1(\Omega,\Psi)$. 
\item[{\rm c)}] $\USCB{X}$ is separable and satisfies the selection principle
\newline
{\rm S}${}_1(\Omega_{\snula},\Psi_{\snula})$. 
\item[{\rm d)}] $\USCB{X}$ is separable and satisfies the selection principle
\newline
{\rm S}${}_1(\tilde \Omega^{\uparrow},\Psi_{\snula})$.
\item[{\rm e)}] $\CP{X}$ satisfies the selection principle {\rm
S}${}_1(\tilde \Omega^{\uparrow},\tilde\Psi)$. 
\item[{\rm f)}] $\CP{X}$ is separable and satisfies the selection principle\\
{\rm S}${}_1(\Omega_{\snula},\Psi_{\snula})$. 
\item[{\rm g)}] $\CP{X}$ is separable and satisfies the selection principle\\
{\rm S}${}_1(\tilde \Omega^{\uparrow},\Psi_{\snula})$.
\end{enumerate}
\end{corollary}
Since $\USCB{X}$ and $\CPB{X}$ are (sequentially) separable
for a separable metrizable space $X$ (see Theorems \ref{th31} and
\ref{th220}), we get the following result.
\begin{corollary}
For $\Psi=\Omega,\Gamma$ and a~separable metrizable space $X$ the following statements are
equivalent:
\begin{enumerate}
\item[{\rm a)}] $\USCB{X}$ satisfies the selection principle {\rm S}${}_1(\tilde \Omega^{\uparrow},\tilde\Psi)$. 
\item[{\rm b)}] $X$ possesses the covering property {\rm S}${}_1(\Omega,\Psi)$. 
\item[{\rm c)}] $\USCB{X}$ satisfies the selection principle {\rm S}${}_1(\Omega_{\snula},\Psi_{\snula})$. 
\item[{\rm d)}] $\USCB{X}$ satisfies the selection principle {\rm S}${}_1(\tilde \Omega^{\uparrow},\Psi_{\snula})$. 
\item[{\rm e)}] $\CP{X}$ satisfies the selection principle {\rm S}${}_1(\tilde \Omega^{\uparrow},\tilde\Psi)$. 
\item[{\rm f)}] $\CP{X}$ satisfies the selection principle {\rm S}${}_1(\Omega_{\snula},\Psi_{\snula})$. 
\item[{\rm g)}] $\CP{X}$ satisfies the selection principle {\rm S}${}_1(\tilde \Omega^{\uparrow},\Psi_{\snula})$.
\end{enumerate}
\end{corollary}
\par
Let us recall that a~set $F\subseteq X$ is a~{\emm zero set} if there exists a continuous function $f\in\CP{X}$ such that $F=\{x\in X:f(x)=0\}$. A~complement of a~zero set is a~{\emm cozero set}.
We denote by $\Phi_{\scriptstyle cz}(X)$ the family of all $\varphi$-covers of $X$ consisting of cozero sets. A~$\varphi$-cover~${\mathcal U}$ is {\emm functionally shrinkable} if there exists a~$\varphi$-cover~${\mathcal W}$ consisting of*/ zero sets such that \eqref{shrink} holds true. The family of all functionally shrinkable $\varphi$-covers of $X$ will be denoted by $\Phi^{\scriptstyle f\!sh}(X)$. We shall be especially interested in the families $\Phi^{\scriptstyle f\!sh}_{\scriptstyle cz}(X)$. 
\par
One can easily see that Theorem~\ref{Th63} may be modified as 
\begin{theorem}
Assume that $\Phi$ is one of the symbols $\Omega$, $\Gamma$ and $\Psi$ is one of the symbols ${\mathcal O}$, $\Omega$, $\Gamma$. Then for any couple $\langle\Phi,\Psi\rangle$ different from $\langle \Omega,{\mathcal O}\rangle$, a~topological space $X$ is an~{\rm S}$_1(\Phi^{\scriptstyle f\!sh}_{\scriptstyle cz},\Psi)$-space if and only if $\CPB{X}$ satisfies the selection principle S${}_1(\Phi_{\snula},\Psi_{\snula})$.
\par
If a~topological space $X$ possesses the property $(\varepsilon)$ then the equivalence holds true also for the couple $\langle\Omega,{\mathcal O}\rangle$.
\end{theorem}
\par
As a~consequence, if we replace in parts b) of Theorems~\ref{PhiDcon}, \ref{PhiScon} and \ref{PhiPcon} the covering property {\rm S}${}_1(\Phi^{\scriptstyle sh},\Psi)$ by the covering property {\rm S}${}_1(\Phi^{\scriptstyle f\!sh}_{\scriptstyle cz},\Psi)$ then we can omit the condition that the topological space $X$ is normal.
\section{Remarks}
The idea to use upper semicontinuous functions for expressing some covering properties of a~topological space arose in \cite{BL1}. That was M.~Sakai~\cite{Sa3} who immediately exploited this idea. The starting point of our study presented in this paper were the results of~\cite{BL2} and the investigation started in~\cite{Osi}.
\par 
Some equivalences of our Theorems are already known. M.~Scheepers~\cite{Sch4} has prove that a~metric separable space $X$ is an~S${}_1(\Omega,\Omega)$-space if and only if $\CP{X}$ satisfies S${}_1({\mathcal D},{\mathcal D})$.
A.V.~Osipov~\cite{Osi} proves several characterizations of topological spaces $X$ with $\CP{X}$ satisfying S${}_1({\mathcal D},{\mathcal S})$, S${}_1({\mathcal S},{\mathcal D})$ or S${}_1({\mathcal S},{\mathcal S})$.
\par
In~\cite{BL2} similar results are presented for~S${}_{\scriptstyle fin}$. Consequently, one can easily see that plenty of our results may be proved also for S${}_{\scriptstyle fin}$.

\end{document}